\documentclass[reqno,12pt,a4paper]{amsart}

\usepackage{mathrsfs}
\usepackage{amssymb}
\usepackage{amsfonts}
\usepackage{latexsym}
\usepackage{amsthm}

\usepackage{graphicx}
\def\lb{\label}

\newcommand{\er}[1]{\textrm{(\ref{#1})}}

\begin{document}


\renewcommand{\theequation}{\arabic{section}.\arabic{equation}}
\theoremstyle{plain}
\newtheorem{theorem}{\bf Theorem}[section]
\newtheorem{lemma}[theorem]{\bf Lemma}
\newtheorem{corollary}[theorem]{\bf Corollary}
\newtheorem{proposition}[theorem]{\bf Proposition}
\newtheorem{definition}[theorem]{\bf Definition}
\newtheorem{remark}[theorem]{\it Remark}

\def\a{\alpha}  \def\cA{{\mathcal A}}     \def\bA{{\bf A}}  \def\mA{{\mathscr A}}
\def\b{\beta}   \def\cB{{\mathcal B}}     \def\bB{{\bf B}}  \def\mB{{\mathscr B}}
\def\g{\gamma}  \def\cC{{\mathcal C}}     \def\bC{{\bf C}}  \def\mC{{\mathscr C}}
\def\G{\Gamma}  \def\cD{{\mathcal D}}     \def\bD{{\bf D}}  \def\mD{{\mathscr D}}
\def\d{\delta}  \def\cE{{\mathcal E}}     \def\bE{{\bf E}}  \def\mE{{\mathscr E}}
\def\D{\Delta}  \def\cF{{\mathcal F}}     \def\bF{{\bf F}}  \def\mF{{\mathscr F}}
\def\c{\chi}    \def\cG{{\mathcal G}}     \def\bG{{\bf G}}  \def\mG{{\mathscr G}}
\def\z{\zeta}   \def\cH{{\mathcal H}}     \def\bH{{\bf H}}  \def\mH{{\mathscr H}}
\def\e{\eta}    \def\cI{{\mathcal I}}     \def\bI{{\bf I}}  \def\mI{{\mathscr I}}
\def\p{\psi}    \def\cJ{{\mathcal J}}     \def\bJ{{\bf J}}  \def\mJ{{\mathscr J}}
\def\vT{\Theta} \def\cK{{\mathcal K}}     \def\bK{{\bf K}}  \def\mK{{\mathscr K}}
\def\k{\kappa}  \def\cL{{\mathcal L}}     \def\bL{{\bf L}}  \def\mL{{\mathscr L}}
\def\l{\lambda} \def\cM{{\mathcal M}}     \def\bM{{\bf M}}  \def\mM{{\mathscr M}}
\def\L{\Lambda} \def\cN{{\mathcal N}}     \def\bN{{\bf N}}  \def\mN{{\mathscr N}}
\def\m{\mu}     \def\cO{{\mathcal O}}     \def\bO{{\bf O}}  \def\mO{{\mathscr O}}
\def\n{\nu}     \def\cP{{\mathcal P}}     \def\bP{{\bf P}}  \def\mP{{\mathscr P}}
\def\r{\rho}    \def\cQ{{\mathcal Q}}     \def\bQ{{\bf Q}}  \def\mQ{{\mathscr Q}}
\def\s{\sigma}  \def\cR{{\mathcal R}}     \def\bR{{\bf R}}  \def\mR{{\mathscr R}}
\def\S{\Sigma}  \def\cS{{\mathcal S}}     \def\bS{{\bf S}}  \def\mS{{\mathscr S}}
\def\t{\tau}    \def\cT{{\mathcal T}}     \def\bT{{\bf T}}  \def\mT{{\mathscr T}}
\def\f{\phi}    \def\cU{{\mathcal U}}     \def\bU{{\bf U}}  \def\mU{{\mathscr U}}
\def\F{\Phi}    \def\cV{{\mathcal V}}     \def\bV{{\bf V}}  \def\mV{{\mathscr V}}
\def\P{\Psi}    \def\cW{{\mathcal W}}     \def\bW{{\bf W}}  \def\mW{{\mathscr W}}
\def\o{\omega}  \def\cX{{\mathcal X}}     \def\bX{{\bf X}}  \def\mX{{\mathscr X}}
\def\x{\xi}     \def\cY{{\mathcal Y}}     \def\bY{{\bf Y}}  \def\mY{{\mathscr Y}}
\def\X{\Xi}     \def\cZ{{\mathcal Z}}     \def\bZ{{\bf Z}}  \def\mZ{{\mathscr Z}}
\def\O{\Omega}

\newcommand{\gA}{\mathfrak{A}}
\newcommand{\gB}{\mathfrak{B}}
\newcommand{\gC}{\mathfrak{C}}
\newcommand{\gD}{\mathfrak{D}}
\newcommand{\gE}{\mathfrak{E}}
\newcommand{\gF}{\mathfrak{F}}
\newcommand{\gG}{\mathfrak{G}}
\newcommand{\gH}{\mathfrak{H}}
\newcommand{\gI}{\mathfrak{I}}
\newcommand{\gJ}{\mathfrak{J}}
\newcommand{\gK}{\mathfrak{K}}
\newcommand{\gL}{\mathfrak{L}}
\newcommand{\gM}{\mathfrak{M}}
\newcommand{\gN}{\mathfrak{N}}
\newcommand{\gO}{\mathfrak{O}}
\newcommand{\gP}{\mathfrak{P}}
\newcommand{\gQ}{\mathfrak{Q}}
\newcommand{\gR}{\mathfrak{R}}
\newcommand{\gS}{\mathfrak{S}}
\newcommand{\gT}{\mathfrak{T}}
\newcommand{\gU}{\mathfrak{U}}
\newcommand{\gV}{\mathfrak{V}}
\newcommand{\gW}{\mathfrak{W}}
\newcommand{\gX}{\mathfrak{X}}
\newcommand{\gY}{\mathfrak{Y}}
\newcommand{\gZ}{\mathfrak{Z}}

\def\ve{\varepsilon}   \def\vt{\vartheta}    \def\vp{\varphi}    \def\vk{\varkappa}

\def\Z{{\mathbb Z}}    \def\R{{\mathbb R}}   \def\C{{\mathbb C}}    \def\K{{\mathbb K}}
\def\T{{\mathbb T}}    \def\N{{\mathbb N}}   \def\dD{{\mathbb D}}


\def\la{\leftarrow}              \def\ra{\rightarrow}            \def\Ra{\Rightarrow}
\def\ua{\uparrow}                \def\da{\downarrow}
\def\lra{\leftrightarrow}        \def\Lra{\Leftrightarrow}


\def\lt{\biggl}                  \def\rt{\biggr}
\def\ol{\overline}               \def\wt{\widetilde}
\def\no{\noindent}


\let\ge\geqslant                 \let\le\leqslant
\def\lan{\langle}                \def\ran{\rangle}
\def\/{\over}                    \def\iy{\infty}
\def\sm{\setminus}               \def\es{\emptyset}
\def\ss{\subset}                 \def\ts{\times}
\def\pa{\partial}                \def\os{\oplus}
\def\om{\ominus}                 \def\ev{\equiv}
\def\iint{\int\!\!\!\int}        \def\iintt{\mathop{\int\!\!\int\!\!\dots\!\!\int}\limits}
\def\el2{\ell^{\,2}}             \def\1{1\!\!1}
\def\sh{\sharp}
\def\wh{\widehat}
\def\bs{\backslash}

\def\all{\mathop{\mathrm{all}}\nolimits}
\def\Area{\mathop{\mathrm{Area}}\nolimits}
\def\arg{\mathop{\mathrm{arg}}\nolimits}
\def\const{\mathop{\mathrm{const}}\nolimits}
\def\det{\mathop{\mathrm{det}}\nolimits}
\def\diag{\mathop{\mathrm{diag}}\nolimits}
\def\diam{\mathop{\mathrm{diam}}\nolimits}
\def\dim{\mathop{\mathrm{dim}}\nolimits}
\def\dist{\mathop{\mathrm{dist}}\nolimits}
\def\Im{\mathop{\mathrm{Im}}\nolimits}
\def\Iso{\mathop{\mathrm{Iso}}\nolimits}
\def\Ker{\mathop{\mathrm{Ker}}\nolimits}
\def\Lip{\mathop{\mathrm{Lip}}\nolimits}
\def\rank{\mathop{\mathrm{rank}}\limits}
\def\Ran{\mathop{\mathrm{Ran}}\nolimits}
\def\Re{\mathop{\mathrm{Re}}\nolimits}
\def\Res{\mathop{\mathrm{Res}}\nolimits}
\def\res{\mathop{\mathrm{res}}\limits}
\def\sign{\mathop{\mathrm{sign}}\nolimits}
\def\span{\mathop{\mathrm{span}}\nolimits}
\def\supp{\mathop{\mathrm{supp}}\nolimits}
\def\Tr{\mathop{\mathrm{Tr}}\nolimits}
\def\BBox{\hspace{1mm}\vrule height6pt width5.5pt depth0pt \hspace{6pt}}
\def\where{\mathop{\mathrm{where}}\nolimits}
\def\as{\mathop{\mathrm{as}}\nolimits}


\newcommand\nh[2]{\widehat{#1}\vphantom{#1}^{(#2)}}
\def\dia{\diamond}

\def\Oplus{\bigoplus\nolimits}



\def\qqq{\qquad}
\def\qq{\quad}
\let\ge\geqslant
\let\le\leqslant
\let\geq\geqslant
\let\leq\leqslant
\newcommand{\ca}{\begin{cases}}
\newcommand{\ac}{\end{cases}}
\newcommand{\ma}{\begin{pmatrix}}
\newcommand{\am}{\end{pmatrix}}
\renewcommand{\[}{\begin{equation}}
\renewcommand{\]}{\end{equation}}
\def\eq{\begin{equation}}
\def\qe{\end{equation}}
\def\[{\begin{equation}}
\def\bu{\bullet}

\title[{Estimates of length of spectrum}]
        {On the measure of the spectrum of direct integrals}
\date{\today}

\def\Wr{\mathop{\rm Wr}\nolimits}
\def\BBox{\hspace{1mm}\vrule height6pt width5.5pt depth0pt \hspace{6pt}}

\def\Diag{\mathop{\rm Diag}\nolimits}

\date{\today}
\author[Anton Kutsenko]{Anton Kutsenko}
\address{Laboratoire de M\'ecanique Physique, UMR CNRS 5469,
Universit\'e Bordeaux 1, Talence 33405, France,  \qqq email \
kucenkoa@rambler.ru }


\maketitle

\begin{abstract}
We obtain the estimate of the Lebesgue measure of the spectrum for
 the direct integral of matrix-valued functions.
 These estimates are applicable for a wide class of discrete periodic
 operators. For example: these results give new and sharp spectral bounds for 1D periodic Jacobi
matrices and 2D discrete periodic Schrodinger operators.
\end{abstract}

\section{Introduction}

Let $K$ be compact metric space with regular Borel measure $\m$,
which satisfies $0<\m(G)<+\iy$ for any open subset $G\ss K$. Define
the Hilbert space $L_N^2=\os_{n=1}^NL^2(K)$ of all
quadratic-summable vector functions $f:K\to\C^N$. Let $\cA: L^2_N\to
L^2_N$ be some bounded self-adjoint operator given by
\[\lb{001}
 \cA(f)(k)=A(k)f(k),\ \ \forall k\in K,
\]
where $A:K\to\C^{N\ts N}$ is some continuous matrix-function. The
matrices $A(k)$ are self-adjoint with dimension $N\ts N$ for any
$k\in K$. It is well known (see e.g. \cite{RS}) that the spectrum of
the operator $\cA$ consists of spectra of $A(k)$
\[\lb{002}
 \s(\cA)=\bigcup_{k\in K}\s(A(k)).
\]
Denoting eigenvalues of $A(k)$ as $\l_1(k)\le...\le\l_N(k)$ we
obtain
\[\lb{003}
 \s(\cA)=\bigcup_{n=1}^N\l_n(K).
\]
Note that all $\l_n:K\to\R$ are continuous functions, since $A$ is a
continuous function on $K$. Our goal is to obtain estimates of the
Lebesgue measure of the spectrum ${\rm mes}(\s(\cA))$ in terms of
$A(k)$, but without calculating eigenvalues $\l_n(k)$. Simple
estimate immediately gives us
\[\lb{004}
 {\rm mes}(\s(\cA))\le 2\|\cA\|=2\max_{k\in K}\|A(k)\|,
\]
but usually this estimate is not very accurate. For example if
$A(k)=A_0=\const$, $k\in K$ then \er{004} gives us $ {\rm
mes}(\s(\cA))\le2\|A_0\|$ but in fact $ {\rm mes}(\s(\cA))=0$.

We restrict our consideration to the case in which the  matrices
$A(k)$ are of the form
\[\lb{005}
 A(k)=\sum_{j=1}^M (\vp_j(k)A_j+\ol{\vp_j(k)}A_j^*),
\]
where $^*$ denotes hermitian conjugate, $\ol{z}$ denotes complex
conjugate, $A_m$ are constant matrices (not necessarily
self-adjoint) and $\vp_m:K\to\C$ are some continuous functions. For
any subset $S\ss\C$ we denote
\[\lb{006}
 \diam(S)=\sup\{|z_1-z_2|,\ z_1,z_2\in S\}.
\]

\begin{theorem}\lb{T1}
For the operator $\cA$ \er{001} with $A(k)$ satisfied \er{005} the
following estimate for the Lebesgue measure of the spectrum is
fulfilled
\[\lb{007}
 {\rm
 mes}(\s(\cA))\le2\sum_{j=1}^M\diam(\vp_j(K))\Tr(A_j^*A_j)^{\frac12}.
\]
\end{theorem}

This Theorem can be extended to the case of "piecewise" continuous
functions $\vp_m$. The bound \er{007} can be improved when some
terms in \er{005} have overlapping spectra. To find absolute gap in
the spectrum of $\cA$ we can apply estimates of spectral curves
$\l_n^-\le\l_n(k)\le\l_n^+$ \er{109}, where $\l_n^{\pm}$ does not
depend on $k$.

Usually, the discrete Schr\"odinger operators on periodic graphs are
unitarily equivalent to the direct integrals of matrices $A(k)$ of
the form \er{005}. This allows us to obtain efficient bounds of the
total length of the spectral bands of such operators. Now we
consider some of the most common discrete periodic operators.

{\bf  1D Jacobi matrices with matrix valued periodic coefficients.}
Consider $1D$ periodic Jacobi matrix
$\cJ:\ell^2_m(\Z)\to\ell^2_m(\Z)$ given by
\[\lb{201}
 \cJ y={a}^*_{n-1}y_{n-1}+b_ny_n+a_ny_{n+1},\ \
 y_n\in\C^m,
\]
where matrices $a_n,b_n\in\C^{m\ts m}$ satisfy
$$
 a_{n+p}=a_n,\ \ b_{n+p}=b_n,\ \ b_n=b_n^*,\ \ {\rm for\ all}\ \ n\in\Z
$$
for some period $p\in\N$. This operator is unitarily equivalent to
$\cA$ \er{202}. Applying Theorem \ref{T1} we obtain
\begin{theorem}\lb{T2}
The following estimate is fulfilled
\[\lb{205a}
 {\rm mes}(\s(\cJ))\le4\min_n\Tr(a_n^*a_n)^{\frac12}.
\]
\end{theorem}
For $m=1$ \er{205a} improves the well known estimate from \cite{DS}
(and \cite{PR} for the almost periodic case)
\[\lb{206}
 {\rm mes}(\s(\cJ))\le4|a_1...a_p|^{\frac1p}.
\]
For example if $p>1$ and all $a_n=T$  except one $a_{n_0}=1$ inside
the period, then \er{205a} gives ${\rm mes}(\s(\cJ))\le1$ while
\er{206} gives ${\rm mes}(\s(\cJ))\le T^{\frac{p-1}p}\to\iy$ for
$T\to\iy$.

For $m>1$ we are unaware of existence in the literature of any
estimates similar to \er{205a}. Moreover, estimate \er{205a} is
sharp. In particular, inequality \er{205} becomes the equality in
the case of $a_n=I_m$ (identity matrix) and $b_n=\diag(4k)_1^m$ for
all $n$ (see example below \er{205}).

{\bf  2D discrete periodic Schrodinger operator.} Consider the
operator
\[\lb{s001}
 \cJ_1 y_{n,m}=y_{n,m-1}+y_{n,m+1}+y_{n-1,m}+y_{n+1,m}+q_{n,m}y_{n,m},\
 \ \ n,m\in\Z
\]
with $N,M$-periodic sequence $q_{n,m}$, i.e.
\[\lb{s002}
 q_{n+N,m+M}=q_{n,m}\in\R\ \ {\rm for\ all}\ n,m\in\Z.
\]
The operator $\cJ_1$ is unitarily equivalent to the multiplication
by the matrix $J_1(k_1,k_2)$ \er{jk1k2}. Applying Theorem \ref{T1}
to this case leads to
\begin{theorem}\lb{T3}
For any $N,M$-periodic sequence $q_{n,m}\in\R$ the following
estimate is fulfilled
\[\lb{est3}
 {\rm mes}(\s(\cJ_1))\le4(N+M).
\]
\end{theorem}
For $N=M=1$ \er{est3} reaches equality in the case of
$q_{n,m}=\const$.

For arbitrary periods $N,M$ we construct the example \er{pot} with
$|\s(\cJ_1)|\approx4\max(N,M)$, which is 2 times worse than
\er{est3}. This example shows us that the power of $N,M$ in
\er{est3} is precisely $1$ but common factor equal to $4$ in
\er{est3} may probably be reduced. For the moment, it can only be
claimed that the exact value of this common factor lies between $2$
and $4$.

\section{Proof of Theorem \ref{T1}}

For self-adjoint matrices $B$, $C$ we will write $B\le C$ iff $C-B$
is a positive-semidefinite matrix, i.e. $x^*(C-B)x\ge0$ for all
$x\in\C^N$. Also we denote $|B|=(B^*B)^{\frac12}$. The matrix $|B|$
positive-semidefinite and self-adjoint.

\begin{lemma}\lb{L1}
For any complex matrix $B$ the following inequalities are fulfilled
\[\lb{101}
 -|B|-|B^*|\le B+B^*\le|B|+|B^*|.
\]
\end{lemma}
{\it Proof.} There exists unitary matrix $U$ ($U^{-1}=U^*$) which
satisfies $B=U|B|$ (polar decomposition). Then $B^*=|B|U^*$ and
$|B^*|=U|B|U^*$. For any $x,y\in\C^N$ we denote $(x,y)_{1}\ev
x^*|B|y$, which is a Hermitian form. Note that $(y,y)_1\ge0$ for any
$y\in\C^N$. Substituting $y=U^*x+x$, $x\in\C^N$ into $(y,y)_1\ge0$
we obtain
\[\lb{102}
 0\le(U^*x+x,U^*x+x)_1=(U^*x,U^*x)_1+(x,x)_1+(U^*x,x)_1+(x,U^*x)_1=
\]
\[\lb{103}
 x^*U|B|U^*x+x^*|B|x+x^*U|B|x+x^*|B|U^*x=x^*(|B^*|+|B|+B+B^*)x,
\]
which gives us the first inequality in \er{101}. Analogously
substituting $y=U^*x-x$, $x\in\C^N$ into $(y,y)_1\ge0$ we obtain the
second inequality in \er{101}. \BBox

{\bf Proof of Theorem \ref{T1}}. There exist points $s_j\in\C$,
$j=1,...,M$ for which
\[\lb{104}
 \frac12\diam(\vp_m(K))=\max\{|s-s_j|,\ s\in\vp_j(K)\}.
\]
Denote
\[\lb{105}
 B_0=\sum_{j=1}^M (s_jA_j+\ol{s_j}A_j^*).
\]
Using \er{101} and \er{104} we deduce that
\[\lb{106}
 A(k)-B_0=\sum_{j=1}^M
 ((\vp_j(k)-s_j)A_j+(\ol{\vp_j(k)-s_j})A_j^*)\le
\]
\[\lb{107}
 \sum_{j=1}^M|\vp_j(k)-s_j|(|A_j|+|A_j^*|)\le\frac12\sum_{j=1}^M\diam(\vp_j(K))(|A_j|+|A_j^*|)\ev B_1
\]
and analogously $-B_1\le A(k)-B_0$. Then we obtain two-sided
inequalities
\[\lb{108}
 B_0-B_1\le A(k)\le B_0+B_1,
\]
where $B_0$, $B_1$ do not depend on $k$. Thus for eigenvalues of
$A(k)$ (see above \er{003}) we deduce that
\[\lb{109}
 \l_n^-\le\l_n(k)\le\l_n^+,
\]
where $\l_1^{\pm}\le...\le\l_N^{\pm}$ are eigenvalues of $B_0\pm
B_1$ respectively. Identity \er{003} gives us
\[\lb{110}
 \s(\cA)=\bigcup_{n=1}^N\l_n(K)\ss\bigcup_{n=1}^N[\l_n^-,\l_n^+].
\]
Since $\l_n^{\pm}$ does not depend on $k$, we obtain
\[\lb{111}
 {\rm
 mes}(\s(\cA))\le\sum_{n=1}^N(\l_n^+-\l_n^-)=\Tr(B_0+B_1)-\Tr(B_0-B_1)=2\Tr
 B_1.
\]
Combining \er{111} with definition of $B_1$ \er{107} and with
identity $\Tr|C|=\Tr|C^*|$ for any complex matrix $C$ gives us
\er{007}. \BBox

\section{Proof of Theorem \ref{T2}}

The operator $\cJ$ \er{201} is unitarily equivalent to the operator
$\cA:L^2_{mp}\to L^2_{mp}$ (see e.g. \cite{KKu}, \cite{KKu1}), where
$L_{mp}^2=\os_{n=1}^{mp}L^2[0,2\pi]$ and
\[\lb{202}
 \cA(f)(k)=A(k)f(k),\ \ f\in L^2_{mN},\ \ A(k)=\left(\begin{array}{ccccc} b_1 & a_1 & 0 & ... & e^{-ik}a_0^* \\
                                      a_1^* & b_2 & a_2 & ... & 0 \\
                                      0 & a_2^* & b_3 & ... & 0 \\
                                      ... & ... & ... & ... & ... \\
                                      e^{ik} a_0 & 0 & 0 & ... &
                                      b_{N}
  \end{array}\right).
\]
Then $A(k)=A_0+\vp_0(k)B_0+\ol{\vp}_0(k)B^*_0$, where $\vp_0=e^{ik}$
and constant matrices $A_0$, $B_0$ are given by
\[\lb{203}
 A_0=\left(\begin{array}{ccccc} b_1 & a_1 & 0 & ... & 0 \\
                                      a_1^* & b_2 & a_2 & ... & 0 \\
                                      0 & a_2^* & b_3 & ... & 0 \\
                                      ... & ... & ... & ... & ... \\
                                      0 & 0 & 0 & ... &
                                      b_{N}
  \end{array}\right),\ \ B_0=\left(\begin{array}{ccccc} 0 & 0 & 0 & ... & 0 \\
                                      0 & 0 & 0 & ... & 0 \\
                                      0 & 0 & 0 & ... & 0 \\
                                      ... & ... & ... & ... & ... \\
                                      a_0 & 0 & 0 & ... &
                                      0
  \end{array}\right).
\]
Applying Theorem \ref{T1} to this case we deduce that
\[\lb{204}
 {\rm mes}(\s(\cA))\le2{\rm
 diam}(\vp_0([0,2\pi]))\Tr(B_0^*B_0)^{\frac12}=4\Tr(a_0^*a_0)^{\frac12}.
\]
Since we may shift sequences $a_n$, $b_n$, it follows that any
element $a_n$ may be chosen instead of $a_0$ and thus we obtain
\[\lb{205}
 {\rm mes}(\s(\cJ))={\rm
 mes}(\s(\cA))\le4\min_n\Tr(a_n^*a_n)^{\frac12},
\]
which coincides with \er{205a} in Theorem \ref{T2}. Estimate
\er{205} is sharp. Let $J$ be Jacobi matrix with elements $a_n=I_m$
($m\ts m$ identity matrix) and $b_n=\diag(4k)_{k=1}^m$ for any $n$.
Since all $a_n$ and $b_n$ are diagonal matrices, $J$ is unitarily
equivalent to the direct sum of scalar Jacobi operators. In our case
this is the direct sum of shifted discrete Shr\"odinger operators
$\os_{k=1}^m(J^{0}+4kI)$ ($J_0$ is a scalar Jacobi matrix with
$a_n^0=1$, $b_n^0=0$ and $I$ is identity operator). Then
$$
 \s(
 J)=\bigcup_{k=1}^{m}\s(J^{0}+4kI)=\bigcup_{k=1}^{m}[-2+4k,2+4k]=[2,2+4m],
$$
which gives us ${\rm mes}(\s(J))=4m=4\Tr (a_na_n^*)^{\frac12}$.

{\bf Remark.} Now we restrict the operator $\cA$ on some interval
$[\a,\b]\ss[0,\pi]$, i.e. consider $\cA_{\a\b}:
\os_{n=1}^{mN}L^2[\a,\b]\to\os_{n=1}^{mN} L^2[\a,\b]$ given by
\er{202}. Then Theorem \ref{T1} gives us
\[\lb{207}
 {\rm
 mes}(\s(\cA_{\a\b}))\le4\sin\frac{\b-\a}2\min_n\Tr(a_n^*a_n)^{\frac12},
\]
since $\diam(\vp_0([\a,\b]))=2\sin\frac{\b-\a}2$.

\section{Proof of the Theorem \ref{T3}}

Introduce the following self-adjoint $N\ts N$ matrices
\[\lb{s1}
 S(k_1)=\ma 0 & 0 & ... & e^{-2\pi ik_1} \\
            0 & 0 & ... & 0 \\
            ... & ... & ... & ...  \\
            e^{2\pi i k_1} & 0 & ... & 0
             \am,\ \
    A_m=\ma q_{1,m} & 1 & ... & 0 \\
            1 & q_{2,m} & ... & 0 \\
            ... & ... & ... & ...  \\
            0 & 0 & ... & q_{N,m}
             \am,\ \ k_1\in\R,
\]
here the matrix $S$ contains two non-zero components and $A$ is a
tridiagonal matrix. Also introduce the following self-adjoint $NM\ts
NM$ matrices
\[\lb{r1b}
 R(k_2)=\ma 0 & 0 & ... & e^{-2\pi ik_2}I \\
            0 & 0 & ... & 0 \\
            ... & ... & ... & ...  \\
            e^{2\pi i k_2}I & 0 & ... & 0
             \am,\ \ B=\ma A_1 & I & ... & 0 \\
            I & A_2 & ... & 0 \\
            ... & ... & ... & ...  \\
            0 & 0 & ... & A_M
             \am,\ \ k_2\in\R,
\]
\[\lb{jk1k2}
 C(k_1)=\diag(S(k_1))_{1}^M,\ \ J_1(k_1,k_2)=B+C(k_1)+R(k_2).
\]
We apply the standard scheme of rewriting periodic operator into the
direct integral of operators with discrete spectrum (see e.g.
\cite{RS}, XIII.16, p.279). Define the following unitary operator
\[\lb{unit}
 \cU:\ell^2(\Z^2)\to\int\limits_{k_1,k_2\in[0,1]}^{\os}\C^{NM}dk_1dk_2,\
 \ \cU(u_{n,m})=(v_{n,m}(k_1,k_2))_{1,1}^{N,M},
\]
where
\[\lb{unit1}
 v_{n,m}(k_1,k_2)=\sum_{n_1,m_1\in\Z}\exp(2\pi
 i(n_1k_1+m_1k_2))u_{n+n_1N,m+m_1M}.
\]
It can be shown that the operator
$\cJ_1:\ell^2(\Z^2)\to\ell^2(\Z^2)$ \er{001} is unitarily equivalent
to the operator of multiplying by the matrix $J(k_1,k_2)$
\[\lb{unit2}
 \cU\cJ_1\cU^{-1}v=J_1(k_1,k_2)v,\ \ \forall v=(v_{n,m}(k_1,k_2))_{1,1}^{N,M}.
\]
By the analogy with the Proof of Theorem \ref{T2}, applying the
Theorem \ref{T1} to the operator of multiplication by the matrix
$J_1(k_1,k_2)$ (see \er{s1}-\er{jk1k2}) leads to the results of
Theorem \ref{T3}.

\section{2D discrete periodic Schrodinger operator with large spectrum}

Without loss of generality we assume $M\ge N$. Introduce the
potential
\[\lb{pot}
 q_{n,m}=\ve^{-1}m,\ \ n=1,...,N,\ \ m=1,...,M,\ \ \ve>0.
\]
Then the matrix $J_1$ \er{jk1k2} is
\[\lb{J1}
 J_1(k_1,k_2)=\ve^{-1}(J_0+\ve D(k_1,k_2)),\ \ {\rm where}\ \
 J_0=\diag(mI)_1^M,
\]
\[\lb{J2}
 D(k_1,k_2)=\ma S_1(k_1) & I & ... & e^{-2\pi ik_2} \\
            I & S_1(k_1) & ... & 0 \\
            ... & ... & ... & ...  \\
            e^{2\pi ik_2} & 0 & ... & S_1(k_1)
             \am,\ \ S_1=\ma 0 & 1 & ... & e^{-2\pi ik_1} \\
            1 & 0 & ... & 0 \\
            ... & ... & ... & ...  \\
            e^{2\pi ik_1} & 0 & ... & 0
             \am.
\]
For identifying the spectrum of $J$ we will apply the regular
perturbation theory with small parameter $\ve$. The eigenvalues of
$J_0$ \er{J1} are
$$
 \l_{n}^{(j)}(0)=j,\ \ n=1,...,N,\ \ j=1,...,M,
$$
with the corresponding eigenvectors
$$
 e_n^{(j)}=e_{n+N(j-1)},\ \ {\rm where}\ \ e_i=(\d_{ij})_1^{NM}
$$
and $\d$ is a Kronecker symbol. Each eigenvalue $j$ has multiplicity
$N$. The perturbation theory for multiple eigenvalues (see e.g.
\cite{K}) tells us that the eigenvalues of $J$ \er{J1} satisfy
\[\lb{peig1}
 \l_{n}^{(j)}(\ve)=\ve^{-1}(j+\ve\wt\l_{n}^{(j)}+O(\ve^2)),
\]
where $\wt\l_{n}^{(j)}\ev\wt\l_{n}^{(j)}(k_1,k_2)$, $n=1,...,N$ are
eigenvalues of the matrix
$$
 (e_1^{(j)}..e_N^{(j)})^{\top}D(k_1,k_2)(e_1^{(j)}..e_N^{(j)})=S_1.
$$
It is well known that the spectrum of $S_1(k_1)$, $k_1\in[0,1]$
coincides with $[-2,2]$ (since $S_1$ corresponds to the 1D discrete
non-perturbed Schrodinger operator), i.e.
\[\lb{pspec}
 \bigcup_{n=1}^N\bigcup\limits_{k_1,k_2\in[0,1]}\{\wt\l_{n}^{(j)}(k_1,k_2)\}=[-2,2].
\]
So, using \er{peig1} with \er{pspec} we deduce that
$$
 \bigcup_{n,j=1}^{N,M}\bigcup\limits_{k_1,k_2\in[0,1]}\{\l_{n}^{(j)}(\ve)\}=4M+O(\ve),
$$
since intervals $[-2,2]+\ve^{-1}j$ do not overlap each other for
sufficiently small $\ve$. Then the spectrum of $\s(\cJ_1)$ for
$\cJ_1$ with potential \er{pot} has Lebesgue measure
$$
 {\rm mes}(\s(\cJ_1))=4M+O(\ve)=4\max\{N,M\}+O(\ve).
$$

{\bf Acknowledgements.} I would to express gratitude to Prof. E.
Korotyaev for stimulating discussions and helpful remarks, to Prof.
B. Simon for useful comments and for the reference to the paper
\cite{PR}, and to Prof. Michael J. Gruber, who has advised me to use
$\Tr(a_na_n^*)^{\frac12}$  instead of $\|a_n\|\rank a_n$ in
\er{205a} (which was in the first version of this paper).

\end{document}